\documentclass[11pt, a4paper, DIV11, oneside, final]{amsart}

\usepackage{verbatim}
\usepackage[latin1]{inputenc}
\usepackage[T1]{fontenc}
\usepackage[stretch=10]{microtype}
\usepackage{tikz}
\usetikzlibrary{calc}
\usetikzlibrary{intersections}
\usepackage{hyphenat}
\usepackage{fixltx2e}
\usepackage{wrapfig}
\usepackage{mathtools}
\usepackage{hyphenat}

\usepackage{amssymb}
\usepackage{tabularx}

\usepackage[DIV11, headinclude=true]{typearea}
\usepackage[hidelinks, linktoc=all]{hyperref}

\newcommand{\ovB}{\bar{B}}
\newcommand{\BB}{\mathcal{B}}
\newcommand{\CC}{\mathcal{C}}
\newcommand{\DD}{\mathcal{D}}

\newcommand{\FFF}{\mathbf{F}}

\newcommand{\LL}{\mathcal{L}}

\newcommand{\UU}{\mathcal{U}}
\newcommand{\WW}{\mathcal{W}}

\newcommand{\id}{{\rm id}}

\newcommand{\divv}{{\rm div}}
\newcommand{\gradd}{{\rm grad}}
\newcommand{\Det}{{\rm Det}}
\newcommand{\real}{\mathbb{R}}
\newcommand{\MMM}{\mathbb{M}}
\newcommand{\complex}{\mathbb{C}}
\newcommand{\integer}{\mathbb{Z}}

\newcommand{\N}{\mathbb{N}}

\newcommand{\VV}{\mathcal{V}}
\renewcommand{\tilde}{\widetilde}

\newcommand{\cone}{\mathbf{C}}
\newcommand{\sphere}{\mathbf{S}^{d-1}}

\DeclareMathOperator{\supp}{supp} \DeclareMathOperator{\Id}{Id}

\newtheorem{proposition}{Proposition}[section]

\newtheorem{theorem}[proposition]{Theorem}

\newtheorem{definition}[proposition]{Definition}

\theoremstyle{remark}
\newtheorem{remark}[proposition]{Remark}
\theoremstyle{definition}

\numberwithin{equation}{section}

\begin{document}
\author{Viviane Baladi} 
\address{Sorbonne Universit\'e, UPMC Univ. Paris 6, CNRS, Institut de Math\'ematiques de Jussieu\\ 
\phantom{De}(IMJ-PRG), 
4, Place Jussieu, 75005 Paris, France}
\email{viviane.baladi@imj-prg.fr}

\author{Tobias Kuna}
\address{Department of Mathematics and Statistics,
University of Reading, RG66AX, Reading, UK}
\email{t.kuna@reading.ac.uk}

\author{Valerio Lucarini}
\address{Department of Mathematics and Statistics, University of Reading, RG66AX, Reading, UK\\ \phantom{D}
CEN, University of Hamburg, Hamburg 20144, Germany}
\email{v.lucarini@reading.ac.uk}

\title[Response for hyperbolic attractors and discontinuous observables]{Linear and fractional response for the SRB measure of smooth hyperbolic attractors and discontinuous observables}

\begin{abstract}
We consider a smooth one-parameter family
$t\mapsto (f_t:M\to M)$ of diffeomorphisms with compact transitive Axiom A attractors $\Lambda_t$, denoting
by $d \rho_t$ the   SRB measure of $f_t|_{\Lambda_t}$.
Our first result  is that for any function $\theta$ in the  Sobolev space
$H^r_p(M)$, with $1<p<\infty$ and  $0<r<1/p$, the map $t\mapsto \int \theta\, d\rho_t$ is $\alpha$-H\"older continuous for
all $\alpha <r$.  This  applies to
 $\theta(x)=h(x)\Theta(g(x)-a)$ (for all $\alpha <1$) for $h$ and $g$
 smooth and $\Theta$ the
 Heaviside function, if $a$ is not a critical value of $g$. Our second result says that
for any such function  $\theta(x)=h(x)\Theta(g(x)-a)$ so that in addition
the intersection of $\{ x\mid g(x)=a\}$  with the support of $h$ is
 foliated by ``admissible stable leaves'' of $f_t$, the map $t\mapsto \int \theta\, d\rho_t$ 
 is differentiable. (We provide distributional linear response and fluctuation-dissipation formulas for the derivative.)
 Obtaining linear response or fractional response for such observables $\theta$ is motivated
 by extreme-value theory.
\end{abstract}

\date{\today}
\thanks{This work was started in 2013 at the Newton Institute during the programme Mathematics for the
Fluid Earth, and continued in 2014 during  a visit of TK and VL to DMA-ENS, where VB was working at
the time, financed by an Action incitative ENS. VL acknowledges funding from the DFG's Cluster of Excellence CliSAP and from the ERC Starting Investigator Grant NAMASTE - Thermodynamics of the Climate System (Grant agreement No. 257106). VB thanks S. Gou\"ezel and M. Todd for useful comments. We thank the referee for helpful remarks and encouraging us to improve the presentation.}
\maketitle

\section{Introduction}

Response theory has long been an essential ingredient of statistical
physics, because it  provides reliable and intuitive tools to
describe the changes in the physical invariant measure of a 
system resulting, e.g., from the introduction of a forcing or of a parametric modulation,
in terms of some statistical properties of the unperturbed system. The
theory --- heuristic and not mathematically rigorous --- was
initially
proposed by Kubo \cite{K} in a seminal paper on statistical
mechanical systems weakly driven outside of the thermodynamical
equilibrium described by the equilibrium measure in the canonical ensemble. In this
case, it
is possible to use response theory to derive the so-called
fluctuation-dissipation theorem which provides a fundamental
connection between the free and forced (linear) fluctuations of a system \cite{K} and has an enormous importance in a multitude of fields of
physical sciences \cite{MPRV}. Response formulas describe the change in the
expectation value of a given observable $\theta$ of the system.

Response theory has been extremely successful. However
  the  limits involved in constructing the
response formulas were not always proved to be well-defined. Ruelle \cite{RLR,RLRb, RLR9} (see also \cite{J}) established a rigorous response theory for transitive
Axiom A 
\footnote{According to the \textit{chaotic hypothesis} of Gallavotti and Cohen \cite{GC}, Axiom A systems
can be taken as effective models for  dynamical systems with many
degrees of
freedom, and they hence 
have  physical relevance.}  attractors. By definition, Axiom A attractors 
\footnote{The terms Axiom A attractor and (uniformly)
hyperbolic attractor are indeed used interchangeably
in the literature and in the present paper.}
are uniformly hyperbolic, and it is well known that
they possess an invariant SRB (physical) measure. 
The SRB measure is  in general singular with respect to Lebesgue
measure (so, from a physical point of view, the system can be
viewed as out of equilibrium), but 
Ruelle showed that,  {\it when used to average
a smooth enough observable $\theta$,} this measure is nevertheless differentiable with
respect to the parameters controlling small modifications  of the
system, and he provided explicit expressions for the derivatives. 
Ruelle's theory has been successfully applied to predict the response
to perturbations of out of equilibrium systems  \cite{LS, R02, EHL04, CS07}, including  the very high-dimensional case of a climate model \cite{LBHRPW}.
Applying  Ruelle's formulas to  the perturbation of a statistical mechanical system whose statistics is described by the canonical ensemble (in contact with a reservoir at constant temperature), the classical results of Kubo can be straightforwardly recovered. 

Modern spectral methods (see e.g. \cite{GL1}), based on a transfer (Ruelle--Perron--Frobenius)
operator acting on a suitable Banach space of anisotropic distributions,
have provided new proofs
and extensions of Ruelle's results. 
However, up to now, 
all existing rigorous works (in particular \cite{RLR, GL1})
establishing linear response for hyperbolic systems
require smooth (at least $C^1$) observables $\theta$. 
(The situation is very different for locally expanding dynamics, where bounded observables
are more than enough, see e.g. \cite{B14} and references therein. This is due
to the fact that the transfer operator there acts on Banach spaces of
functions and not distributions.)
The goal of the present short note\footnote{In this paper, we only consider
discrete-time dynamical systems, i.e., iterations of a map $f:M\to M$, or of
one-parameter families of maps $f_t:M\to M$.}
   is to investigate response
theory of Axiom~A systems for non differentiable observables, allowing 
in particular discontinuities in the observable. 
Our results below (Theorems~\ref{main} and ~\ref{main2})
apply  to observables involving the Heaviside function $\Theta(x)$, for example
 $\theta(x)=\Theta(g(x))h(x)$, with $g$ and $h$ smooth, assuming  that
 zero is not a critical value of $g$. 
Note that the expectation\footnote{Replacing $\Theta$ by a Dirac mass $\delta_D$ also gives interesting observables, see  Remarks~\ref{basic} and~\ref{blueprint}.}
value of $\Theta(g(x))$ gives the fraction of the total measure where $g$ has
positive value, while the expectation value of $h(x)\Theta(g(x))$
divided by the expectation value of $\Theta(g(x))$ gives the average of $h$
conditional on the fact that $g$ is positive. 
Therefore, such discontinuous observables have a clear probabilistic and physical interpretation. Here is an example
regarding the analysis of extreme events \cite{Letal16} in chaotic dynamical systems: In \cite{LFWK} it is suggested that one can construct estimates of the parameters describing the 
statistics of extremes of a general observable $g$ for Axiom~A systems through suitable combinations of the expectation value of quantities of the form  $(g(x)-a)^n\Theta(g(x)-a)$, for $n\in \N$, where $a$ is a (noncritical)
threshold close enough to the global maximum of $g$ on the attractor of the system if $g$ is smooth, while
 $a$ tends to infinity if $g$ is unbounded. Therefore,  some regularity in the dependence of such discontinuous observables  on perturbations of the dynamics will give control on the response of extreme events --- an obviously relevant practical problem.
 (We refer in particular to \cite[Section III]{LFWK}, see
  Remark~\ref{extrr} below. The discussion around \eqref{quot} below mentions related
  open questions.)

\medskip

In the present work, we consider one-parameter families
of transitive compact hyperbolic attractors $f_t:\Lambda_t\to \Lambda_t$, 
with $\Lambda_t$ contained in a Riemann manifold $M$, denoting the SRB measure
by $\rho_t$, and we obtain two main results:
\smallskip

In Section~\ref{sec1}, we prove Theorem~\ref{main}, which says that for any 
$1<p<\infty$, any $0<r<1/p$, and any observable $\theta$ in the generalised (isotropic, ordinary) 
Sobolev space $H^r_p(M)$, the map 
\begin{equation}
\label{DT}t \mapsto \int \theta d\rho_t
\end{equation}
is $\alpha$-H\"older
continuous for all exponents $\alpha <r$. (Remark~\ref{extrr} shows that this  result
applies to observables of the type $h(x)\Theta(g(x)-a)$, and in that case $r$ can be taken arbitrarily close to $1$.)
We expect that it is possible to find examples of hyperbolic attractors $f_t$ and
observables
$\theta \in H^r_p(M)$ for which  the map \eqref{DT} is not differentiable and maybe to show that $\alpha=r$ is the optimal
regularity limitation.
This is not the first time that a {\it fractional} (i.e., weak) formulation of the response has
been obtained: It is known that, in nonuniformly hyperbolic situations
where  strong bifurcations are present,
differentiability of the response can be violated (even for $C^\infty$ observables).
In fact, $1/2$-H\"older upper {\it and lower} bounds
have been obtained in the setting of the logistic family (i.e., unimodal interval maps, see
 \cite{B14,BBS} and references therein). 
 We emphasize that it is often however extremely difficult to distinguish numerically
 between differentiability (strong, linear, response) and H\"older continuity
 (weak, fractional, response), as is clearly
 discussed in the illuminating preprint
\cite{GWW} of Gottwald, Wormell, and  Wouters, which came to our attention while we were
finishing the present work.
 
\smallskip
In Section~\ref{sec2}, our second main result,
Theorem~\ref{main2}, gives ordinary (i.e., strong) linear response 
for observables $h(x)\Theta(g(x)-a)$, whose singularity set is ``compatible'' with 
the stable cone of the transitive attractor. The precise compatibility condition
involves ``admissible stable leaves,''  (see Definition~\ref{notat}), 
and requires
that the set 
$$\WW_a=\{ x\mid g(x)=a\} \cap \supp(h)$$ 
be
 foliated by admissible stable leaves.
 This condition, which implies transversality of the normal vector to $\WW_a$
 with the ``wave-front set'' of the SRB measure,
  will not surprise experts, but to our knowledge it is the first time that it
 is written down.  By ordinary linear response, we mean here
 differentiability of the map \eqref{DT}, together with  formulas
for the derivative, both of linear response-type \eqref{LRF3} and fluctuation-dissipation-type
 \eqref{FDT}  (both formulas being in the sense of distributions!).
 We explain next how to interpret   (distributional) fluctuation-dissipation
 expressions  
 in hyperbolic settings, if the SRB measure $d\rho_0$ of $f_0$
 is singular with respect to Lebesgue. The starting point is 
 that $d\rho_0$ can be approached (in anisotropic norm)
 by a sequence $\rho^{(n)}_0 d m$ of absolutely continuous measures, 
 with smooth densities, either by
 mollification, or by iterating Lebesgue measure under the dynamics $f^n_0$. Then,
 setting $X_0=\partial_t f_t |_{t=0}\circ f_0^{-1}$, 
our results on the distributional fluctuation-dissipation
 formula \eqref{FDT} easily imply that the terms of the sum
 $$
 -\sum_{k=0}^\infty
 \int (\theta \circ f_0^k ) [\divv (X_0)\rho_0^{(n)}+ \langle X_0 ,\gradd \rho^{(n)}_0\rangle ]\, dm
 $$
decay exponentially (uniformly in $n$), and the sum 
converges to $\partial_t \int \theta \, d\rho_t|_{t=0}$ as $n\to \infty$,
 if the observable $\theta$ lies in the dual of a suitable anisotropic space.
(Admittedly,  $\gradd \rho^{(n)}_0$ is not very nice if $\rho^{(n)}_0$ is constructed by
iterating Lebesgue measure, it is slightly easier to control if $\rho^{(n)}_0$
is obtained by convolution with an approximation of the Dirac mass.)
 
 We would like to emphasize that, in hyperbolic situations (as opposed
to
 locally expanding maps),  anisotropic
norms are still required   to get exponential convergence to equilibrium or
linear response {\it even if the SRB measure is absolutely continuous with respect
to Lebesgue.}  In particular,  fluctuation-dissipation expressions
such as \eqref{FDT} require the observable to have some modulus of continuity
in the stable cones. However, in the case of an absolutely continuous SRB measure, 
the gradient of the
SRB measure will be contained in a space of arbitrarily mild
anisotropic distributions (in the notation below, instead of
$W^{u-1,s-1}_p$ or $\BB^{u-1, |s-1|}$ for $s<0$, one can use $W^{u-1,\sigma}_p$
or $\BB^{u-1,|\sigma|}$ for $\sigma<0$ arbitrarily close to $0$), and
one can then hope to weaken slightly the condition on the observable. 
The toy-model case of smooth enough locally expanding maps is much easier
to tackle, since then
linear response holds for observables in spaces of distributions, so that  one
could obviously replace the Heaviside function by a Dirac mass $\delta_D$, or derivatives of
a Dirac mass.
(To summarise,  there is a hierarchy of difficulties, with locally expanding easiest, singular SRB
measure hardest, and absolutely continuous SRB measure in the middle.)

\smallskip
At a technical level, we work with  anisotropic Banach
norms. For Theorem~\ref{main}, it is more convenient to use a norm introduced in \cite{BTs}.
One of the lessons of the present work is that, for the fractional
response of Theorem~\ref{main}, only a weak, dynamics-independent, condition on the observable is needed,
and a ``cheap mollification trick'' does the job. 
To prove  Theorem~\ref{main2},   the  geometric anisotropic norms
from \cite{GL1}  are more convenient.  The  compatibility
condition on $\WW_a$ given above implies the transversality
needed to check that the observable lies in the dual of the anisotropic space containing
the gradient of the SRB measure,
which is just what the proof requires.

\section{Fractional response for hyperbolic attractors and observables in $H^r_p(M)$}
\label{sec1}

Let $M$ be a smooth $d$-dimensional Riemann manifold.
Recall  that a  nonempty compact set $\Lambda_t\subset M$ is a (uniformly) hyperbolic attractor
for a diffeomorphism $f_t:\Lambda_t \to\Lambda_t$, if
there exists an open neighbourhood $V_t$ of $\Lambda_t$ so that
 $f_t(\overline{V_t})\subset V_t$ and $\cap_{n\ge 0} f_t^n(V_t)=\Lambda_t$,
and if there exist $C>1$ and $\nu<1$ so that  tangent space $T M|_{\Lambda_t}$ 
over $\Lambda_t$ decomposes into
$E^u \oplus E^s$,  with $E^u$ and $E^s$ two $Df_t$-invariant bundles,
of respective dimensions $d_u\ge 1$ and $d_s\ge 1$, 
so that
$$
\|D_x f^n_t|_{E^s_x\to E^s_{f^n_t(x)}}\| \le C \nu^n \, , \,\,\, 
\|D_x f^{-n}_t|_{E^u_x\to E^u_{f^{-n}_t(x)}}\| \le C \nu^n
\, ,  \, \, 
\forall x \in \Lambda_t\, , \, \forall n\ge 1\, .
$$
The attractor is transitive if there exists a point with a dense orbit.
It is well known (see e.g. \cite{KH})
that if $f_t:\Lambda_t\to \Lambda_t$ is a transitive hyperbolic attractor as just defined,
then $f_t$ has a unique SRB (invariant probability) measure on $\Lambda_t$.
In what follows,  there will always be  a compact subset $\Lambda$ of $M$
containing all the hyperbolic attractors $\Lambda_t$ involved, so for all practical purposes we can
assume that $M$ is compact, denoting by $dm$ normalised Lebesgue
measure on $M$, and we can work with finite systems of smooth 
charts and partitions of unity.

For 
real numbers $r$ and $1<p<\infty$, let $H^r_p(M)$ be the Sobolev 
space of functions $\varphi\in L_p(M)$
 so that, denoting by $\Delta$ the Laplacian in $\real^d$
and by $\tilde \varphi_j$ the
function $\varphi$ in the $j$th chart, the $L_p$ norm of each $(\id +\Delta)^{r/2} (\tilde \varphi_j)$ is finite.
If $r>0$, then $H^r_p(M)$ is the closure of  $C^u(M)$ for any $u>r$, for the same norm, see e.g. \cite{Ba}.
If $r=0$ then $H^r_p(M)=L_p(M)$. 
Also, $H^{-r}_p(M)$ is the dual space to $H^r_{p'}(M)$ if $1/p+1/p'=1$, where duality
is given by the scalar product
$ \int_M \varphi_1 \bar \varphi_2 \, dm$. 
Recall that an element of $H^r_p(M)$ is continuous if $r > d/p$, but beware that
we shall only be interested in the case $r<1/p$.

The following theorem is our first technical result on compact transitive
hyperbolic attractors. Its main interest lies in the fact that
the assumption on the observable $\theta$ does not involve the dynamical system. Also, the
proof (which uses a standard mollification
trick) is remarkably simple.

\begin{theorem}[Fractional response]\label{main} 
Fix $\beta \in (0,1)$. Let $t\mapsto f_t$, for $t\in [-\epsilon_0,\epsilon_0]$, be a $C^{2+\beta}$ family
of  $C^3$  diffeomorphisms $f_t$   on a smooth Riemann manifold $M$, 
so that $f_t$ has a transitive compact 
hyperbolic
attractor $\Lambda_t\subset M$. Let $\rho_t$ be
the (unique) SRB measure of $f_t$ on $\Lambda_t$.  Let $\theta :M \to \complex$ be so that
 $\theta \in H^r_p(M)$ for some $1<p<\infty$ and  $0< r <1/p$. Then there exists $\epsilon_1\in (0, \epsilon_0]$ so that  for
any $\alpha <r$ the map
$$
t\mapsto \int_M \theta\,  d\rho_t
$$
is $\alpha$-H\"older continuous on $[-\epsilon_1, \epsilon_1]$. 
\end{theorem}

An important ingredient of the proof of Theorem~\ref{main}
is   the transfer operator $\LL_t$ associated to $f_t$ by setting
\begin{equation}
\label{defLt}
\LL_t (\varphi)(x)=\frac{\varphi (f_t^{-1}(x))}{|\Det Df_t (f_t^{-1}(x))|}\, ,\quad \varphi \in C^2 \, ,
\end{equation}
where $|\Det Df_t|$ is the Jacobian with respect to Lebesgue measure (so that
the dual of the transfer operator leaves Lebesgue measure on $\Lambda$ invariant).
We shall let this operator act on  anisotropic Banach spaces of
distributions on which it is known to have a spectral gap.

\begin{remark}\label{basic} Previous results on linear response (\cite{RLR}--\cite{J}, \cite[\S 8]{GL1})
restrict to observables $\theta$ which are $C^r$ for some
$r\ge 1$.  As pointed out by the referee, 
Demers and Liverani \cite{DL} obtained {\it fractional} (H\"older) response for two-dimensional piecewise hyperbolic
(or smooth hyperbolic) systems and piecewise H\"older test functions.
An inescapable requirement 
is that $\int \theta d\rho_t$ be well defined. Practically, $\theta$ must lie in the dual space to
an anisotropic Banach space containing $\rho_t$. In particular,  we cannot take $\theta$
to be a Dirac mass on $M$, since  $\rho_t$
cannot integrate such distributions in general (see the proofs of Theorems \ref{main} and \ref{main2} for
more information about $\rho_t$). Note however that one can probably show that $\rho_t$ integrates distributions involving the Dirac
distribution under a condition in the spirit of  Theorem~\ref{main2} below (see \cite[Prop.~4.4]{GL1}, 
although this proposition cannot be applied directly since what we need is a statement on the {\it dual} space).  This inescapable requirement does not suffice for full-fledged linear response, since what 
is relevant there is integration against
$\gradd \rho_t$ and not against $\rho_t$ itself, see Section~\ref{sec2} and Remark~\ref{blueprint}.
\end{remark}
\medskip

\begin{remark}[Application to extreme value theory]\label{extrr}
Let $\Theta:\real \to \{0, 
1\}$ be the Heaviside  step function
$$
\Theta(v)=
0 \,  ,\mbox{ for }  v\le 0\, ,\quad
\Theta(v)=1\, ,  \mbox{ for } v>0\, .
$$
An application of  Theorem~\ref{main} to extreme value theory can be obtained by taking 
$$\theta(x)=\theta_a(x)=\Theta(g(x)-a)$$
for a  $C^2$ function $g:M \to \real$ (we assume for simplicity here that
$M$ is compact) and any   threshold
$a\in \real$ which is not a critical value of $g$. The theorem then gives $\alpha$-H\"older continuity of the response for all  $\alpha <1$.  (Indeed, $\Theta(g(x)-a)$ is just the indicator function
of the set $\DD_a=\{x \in M \mid g(x) > a\}$. We claim that this implies
that $\theta(x)\in H^r_p(M)$ for all $1<p<\infty$ and all $0\le r<1/p$, so that
we can take any $\alpha<1$ in Theorem~\ref{main} by letting $p$ tend to $1$ and $r$ tend to $1/p$. If $a$ is larger than   the global maximum of $g$ or
smaller than  the global minimum of $g$, then the claim is obvious since $\DD_a$ 
coincides with the empty set or with
$M$, respectively. Otherwise, recall (see e.g. \cite{Str}) that the indicator function
of a compact interval belongs to $H^r_p(\real)$ for all $0\le r<1/p$ if $1<p<\infty$.  If $a$ is not a critical value of $g$ then the implicit function theorem
gives that the level set $\WW_a=\{x \in M \mid g(x)= a\}$ is
a finite union of smooth codimension-one submanifolds
and it is the boundary of $\DD_a$. Therefore, by Fubini\footnote{See e.g. the proof of \cite[Lemma 23]{BG1} and the references therein to \cite{Str}. We emphasize however that we only need 
here $\theta$ to belong to a suitable (dual) space and not to be a multiplier.}, the indicator function 
of $\DD_a$ belongs to  $H^r_p(M)$, for all $1<p<\infty$ and all
$0\le  r <1/p$, as claimed, since $M$ is compact.)

As a last item of this remark, assume (without loss of generality) that the global maximum
of $g$ is equal to $0$. For $\theta_b(x)=\Theta(g(x)-b)$,  it would be interesting
 to study $\int \theta_b\, d\rho_t$ as a function of $b$, in particular
in view of considering the  regularity of the limit of the  map\footnote{The quotient \eqref{quot} is also interesting for a function $\theta_a=h(x) \Theta(g(x)-a)$ as
in the statement of Theorem~\ref{main2} below, and the fact that $\rho_t$ is equivalent to Lebesgue measure along
unstable leaves renders  plausible that the limit as $a\to 0$ exists, in view of the proof of Theorem~\ref{main2}.}
\begin{equation}\label{quot}
t \mapsto \frac{\int \theta_{as}\, d\rho_t}{\int \theta_a\, d\rho_t}\, , 
\end{equation} 
for fixed $s\in (0,1)$, 
as $a\to 0$ (so that numerator and denominator both tend to zero, since $0$ is the global
maximum).
Indeed,  for fixed $t$,
proving the existence of the limiting quotient
\eqref{quot} as $a \to 0$  and finding an expression for this limit
are key steps to
construct explicit extreme value laws for the observable $g(x)$. (We refer e.g. to \cite[(3.1.7)]{Letal16}
for
 i.i.d. random variables and to \cite[(4.2.6)]{Letal16}
for dynamical systems.)
For  occurrences of similar quotients in the physics literature,
see   \cite[(12)]{LFW} and  \cite[(5)]{LFWK}, where two different kinds of observables are considered   and the limiting expression is derived under the assumption that the limit \eqref{quot} exists in a suitable sense (\cite[Section III]{LFWK} includes a discussion about
varying $t$).
We also refer to the overview in  \cite[Chapter 8, (8.2.4), 
(8.2.13)]{Letal16}.
(The notation in \cite{Letal16,LFW,LFWK} is  $T=a<0$, and $Z+T=as$, so that 
$Z=T(s-1)>0$:
 This additive formulation does not give a useful  limit as $T$ tends to
zero for fixed $Z$, but it is convenient for experiments or simulations.)
The study of the limiting quotient  map \eqref{quot} is beyond the scope of
this paper, and we postpone it to further works. The first test case would be
that of locally expanding maps.
\end{remark}

\begin{proof}[Proof of Theorem~\ref{main}]
Let $V$ be a common open isolating (in fact, attracting)
neighbourhood for all the $f_t$ (up to reducing the value of $\epsilon_0$):
$\Lambda_t=\cap_{n \in \integer} f^n_t(V)=\cap_{n \in \integer_+}f_t^n(V)$, with $f_t(\overline V)$
a strict subset of $V$. We can assume that the closure of $V$ is compact.
The support of each $\rho_t$ is contained in $V$ by construction.

Fix a finite system of $C^\infty$ charts $\psi_i : U_i \to M$, $i=1, \ldots, N$, with open 
bounded domains $U_i\subset \real^d$,
so that $\overline V \subset \cup_{i=1}^M \psi_i(U_i)$, 
and fix a  $C^\infty$ partition of unity on $V$, i.e.,
functions $\phi_j$ with $\sum_{j=1}^N\phi_j(y) =1$ for all $y\in V$, with the support of
$\phi_i$ compactly contained in  $\psi_i(U_i)$. Let $\eta : \real^d \to \real_+$ be a
$C^\infty$ function supported in $[-1,1]^d$, with $\int \eta(x) \, dx =1$. For any small enough $\varepsilon>0$,  put
$\eta_\varepsilon(v)=\varepsilon^{-d} \eta(v/\varepsilon)$, and  set for $x\in M$
\begin{equation}\label{mol}
\theta_\varepsilon(x)=\MMM_\varepsilon(\theta):=\sum_{i=1}^N  \phi_i(x)\bigl (
 \int_{\real^d} \eta_\varepsilon(v) \theta(\psi_i(\psi_i^{-1}(x)-v))\, dv\bigr )
\, .
\end{equation}
In the rest of the proof of the theorem, we ignore the charts, slightly abusing notation. Since the $\psi_i$ and
$\phi_i$ are $C^\infty$, and there are finitely many of them, this does not cause any problems
(we refer to \cite{BTs}, \cite[\S 5]{BLi}, or \cite{Ba} for details).

In \cite{BTs} it is proved that the transfer operator $\LL_t$ defined by
\eqref{defLt} acting on\footnote{We recall the construction of the space in
Appendix~\ref{appA}.} a Banach space
$W_{p'}^{u,s}(f_t,V)=W_{p', \dagger}^{u,s}(f_t,V)$  of anisotropic distributions has essential spectral radius strictly
smaller than $1$, for all real numbers $1<p'<\infty$ and $u-2<s<0<u<2$. In fact, up to further
reducing $\epsilon_0$, we can use the same space $W_{p'}^{u,s}(f_0,V)$ (noted $W^{u,s}_{p'}$
from now on, for simplicity)
for all $\LL_t$ with $|t|< \epsilon_0$.  Recall that 
$\int \LL_t(\varphi)\, dm=\int \varphi\, dm$ so that the dual of $\LL_t$
fixes Lebesgue measure. By standard arguments
(see \cite{GL2} or \cite{Ba}), the spectral radius of $\LL_t$ on $W^{u,s}_{p'}$ is equal to $1$,
and $1$ is a simple eigenvalue of $\LL_t$, for all $t$.
(The spectral properties listed in this paragraph are often referred to
as Perron--Frobenius properties.)

It follows from \cite[App. A]{BTs} (see also \cite{Ba}) that
we have the following bounded inclusions for any $r'>u$ and $\tilde r >|s|$:
$$
H^{r'}_{p'}(M)\subset W^{u,s}_{p'} \subset H^{-\tilde r}_{p'}(M)
\mbox{ and } H^{\tilde r}_{p}(M)\subset (W^{u,s}_{p'})^* \subset H^{- r'}_{p}(M)\, ,
$$
where $1/p+1/p'=1$. (We shall only use the inclusions involving $\tilde r$.)
Also, the partial derivatives of an element of
$W^{u,s}_{p'}$ belong to $W^{u-1,s-1}_{p'}$.

 Denoting by $\tilde \rho_t$ 
 the fixed point of $\LL_t$ 
in   $W^{u,s}_{p'}$,
it is well-known that 
the distribution $\varphi \mapsto  \tilde \rho_t (\varphi)$ 
defined for $\varphi\in H^{\tilde r}_p(M)$ extends to a nonnegative Radon measure, which is
 $f_t$-invariant. Normalising gives a 
probability measure, which is  in fact the unique SRB measure $d\rho_t$ of $f_t$.
Slightly abusing notation, we shall not distinguish between $\tilde \rho_t(\varphi)$ and 
$\int \varphi d\rho_t$, in particular
we drop the tilde from now on.
Note also the fixed point property implies that $\tilde \rho_t$ also belongs to the (smaller) anisotropic space 
$\tilde W^{u,s}_{p'}$ obtained from $W^{u,s}_{p'}$ by
taking a narrower stable cone and a wider unstable cone.

Now, there exists $C$ so that
each component of  $\gradd \rho_t$ has  $\tilde W^{u-1,s-1}_{p'}$-norm bounded
by $C\|\rho_t\|_{\tilde W^{u,s}_{p'}}$.  Since we can take $u \in (1, 1+\beta)$ and
$u-2<s<0$, the transfer operator $\LL_t$ acting on $W^{u-1,s-1}_{p'}$ also has essential spectral radius strictly
smaller than $1$ and enjoys the Perron-Frobenius spectral properties described above.
For any $C^{1+\beta}$ function $\upsilon$ and any $j$,  since $u-1< \beta$
and $\|\partial_j \upsilon\|_{C^\beta}\le \|\upsilon\|_{C^{1+\beta}}$, we have
\begin{equation}\label{Leibb}
\| (\partial_j \upsilon) \rho_t\|_{W^{u-1,s-1}_{p'}}\le \|\upsilon\|_{C^{1+\beta}}
\|\rho_t\|_{\tilde W^{u-1,s-1}_{p'}}\, , \, \, \| \upsilon (\partial_j \rho_t)\|_{W^{u-1,s-1}_{p'}}
\le \|\upsilon\|_{C^\beta}\|\partial_j \rho_t\|_{\tilde W^{u-1,s-1}_{p'}}\, , 
\end{equation}
using the Leibniz-type bound obtained from applying
 the Lasota--Yorke
estimates in \cite{BTs} to the identity map.
For any smooth vector field $Y$,  integrating by parts
and using that the manifold $M$ is boundaryless,
we find 
\begin{equation}\label{boundaryless}
\int_M \divv (Y) \rho_t + \langle \gradd \rho_t, Y \rangle  =0\, ,
\end{equation}
(where $\langle Z, Y \rangle$ denotes the scalar product of vector fields)
in the sense
of distributions. Hence (recalling that Lebesgue measure
is a simple fixed point of each  dual operator $\LL^*_t$), the resolvent
$(1-\LL_t)^{-1}$ is well-defined when acting on $\divv (Y) \rho_t + \langle  \gradd \rho_t, Y\rangle\in W^{u-1,s-1}_{p'}$.
Recalling \eqref{Leibb},  
the image $(1-\LL_t)^{-1}[\divv (Y) \rho_t + \langle  \gradd \rho_t, Y\rangle]$
is a distribution whose $W^{u-1,s-1}_{p'}$-norm
is bounded by a constant (depending only on $\|Y\|_{C^{1+\beta}}$) times the $\tilde W^{u,s}_{p'}$-norm of 
$\rho_t$. 

We may apply the abstract theorem of  \cite[\S 8]{GL1} to the present  spaces
 $W^{u,s}_p$ (although the spaces $W^{u,s}_p$, taken from  \cite{BTs}, are
 not isomorphic to the spaces used in \cite{GL1}):
Mutatis mutandis, we just follow the arguments in  \cite[\S 5.3, App. A.3]{Ba},
where \cite[\S 8]{GL1} was applied to a slightly different anisotropic space
  introduced in  \cite{BT2}.  Translating the statement from \cite[Theorem 2.7]{GL1}
gives that $t \mapsto \rho_t$ is differentiable when the values
are viewed in  $W^{u-2,s-2}_{p'}$. (This would only allow to show H\"older regularity
in Theorem~\ref{main} for $\alpha <1/2$.) However,   we claim that differentiability
holds in the stronger norm $W^{u'-1,s'-1}_{p'}$, for $u-1<u'<u$ and $s-1<s'<s$, arbitrarily close to
$u$ and $s$, respectively. To check this, we follow the argument of \cite[Theorem~2.38]{Ba},
bootstrapping from the fact that for $\delta \in (0,1)$ there exists $C$ so that for all $\varphi$
\begin{equation}\label{missingl}
\| \LL_t \varphi - \LL_0 \varphi\|_{W_{p'} ^{u,s}}
\le C |t|^{\delta} \|\varphi \|_{W^{u+\delta, s+\delta}_{p'}} \,  .
\end{equation}
The bound \eqref{missingl} can be proved by using the mollifiers defined by \eqref{mol}:
First note that
\begin{align*}
\| \LL_t \varphi - \LL_0 \varphi\|_{W_{p'} ^{u,s}}
&\le
\| \LL_t \MMM_\varepsilon(\varphi) - \LL_0 \MMM_\varepsilon(\varphi)\|_{W_{p'} ^{u,s}}\\
&\qquad+\|\LL_t (\MMM_\varepsilon(\varphi) -\varphi)\|_{W_{p'} ^{u,s}}+\|\LL_0 (\MMM_\varepsilon(\varphi) -\varphi)\|_{W_{p'} ^{u,s}}\\
&\le C |t| \|\MMM_\varepsilon(\varphi)\|_{\widehat W_{p'} ^{u+1,s+1}}+
2C \|\MMM_\varepsilon(\varphi) -\varphi\|_{\widehat W_{p'} ^{u,s}}\, , 
\end{align*}
where the  (larger) anisotropic space $\widehat W^{u,s}_{p'}$ 
are obtained from $W^{u,s}_{p'}$ by taking
a suitable wider stable cone and a  suitable narrower unstable cone.
Second  observe that\footnote{To prove both claims,
use $W^{u,s}_{p',\dagger\dagger}(T,V)$ from \cite[App. A]{BTs}
and the results therein,
since $\|(1+\Delta)^{\delta/2} \varphi\|_{W^{u,s}_{p',\dagger\dagger}}
$ is equivalent to $\|\varphi\|_{W^{u+\delta,s+\delta}_{p',\dagger\dagger}}$. 
The change of cones is needed both
to handle the multiplication by $\phi_i$ in the definition of $\MMM_\varepsilon$ and the
change of charts $\psi_i \circ \psi_j^{-1}$.}
$$
\|\MMM_\varepsilon(\varphi)\|_{\widehat W_{p'} ^{u+1,s+1}}
\le C \varepsilon^{\delta-1} \|\varphi\|_{W_{p'} ^{u+\delta,s+\delta}}
\mbox{ 
and }
\|\MMM_\varepsilon(\varphi)-\varphi\|_{\widehat W_{p'} ^{u,s}}
\le C \varepsilon^{\delta} \|\varphi\|_{W_{p'} ^{u+\delta,s+\delta}} \, .
$$
Finally, taking $\varepsilon=|t|$ proves \eqref{missingl}.

\medskip

Next, using that
$$H^{1+r}_{p}(M)\subset (W^{u'-1,s'-1}_{p'})^*
\mbox{ if } 1+r>|s'-1|\, , 
$$ 
the differentiability of  $t \mapsto \rho_t \in W^{u'-1,s'-1}_{p'}$ gives that
 for any $r >0$ (taking $s'<s<0$ close enough to zero), any $1<p'<\infty$, and any function $\varphi \in H^{1+r}_{p}(M)$ , the map
$$t \mapsto \int \varphi \, d \rho_t$$ 
is differentiable, and its derivative at any $|t_0|<\epsilon_0$ is given by the
following  fluctuation-dissipation formula 
\begin{align}
\partial_t (\int_M \varphi \, d \rho_t )|_{t=t_0}
\label{LRF} & =-
\bigl ( (1-\LL_{t_0})^{-1}[\divv (X_{t_0})\rho_{t_0}+ \langle \gradd \rho_{t_0}, X_{t_0} \rangle ]  \bigr) (\varphi)\, , 
\end{align}
where the $C^{1+\beta}$ vector field $X_{t_0}$ is defined by
\begin{equation}\label{X_t}
X_{t_0} = (\partial_t f_t) |_{t=t_0}\circ f_{t_0}^{-1}\, .
\end{equation}
In particular, the derivative $\partial_t (\int_M \varphi \, d \rho_t )|_{t=t_0}$ depends continuously on $t_0$.

Let $1<p<\infty$ and $0<r<1/p$ be so that $\theta\in H^r_p(M)$. Fix   $|s|<r$
small and $\tilde r \in (|s|,r)$ close to $|s|$.
Let $p'$  be so that $1/p+1/p'=1$.

There exists a constant $C$ (see, e.g., the proof of \cite[Lemma 5.3]{BLi}), so that
for any $\varepsilon >0$ 
\begin{equation}\label{willuse}
\|\theta_\varepsilon\|_{H^{1+\tilde r}_p}\le C \|\theta\|_{H^{r}_p} \varepsilon^{r-1-\tilde r} \, .
\end{equation}
To simplify notation, assume that $t_0=0$.
For fixed small $\varepsilon >0$, decomposing $\theta=\theta_\varepsilon+ (\theta-\theta_\varepsilon)$, it suffices
to estimate the terms $|\int \theta_\varepsilon (d\rho_t-d\rho_0)|$ and 
$|\int (\theta_\varepsilon-\theta) d\rho_t|+|\int (\theta_\varepsilon-\theta) d\rho_0|$.

On the one hand, it follows from the linear response formula \eqref{LRF}, the
mean value theorem, and \eqref{willuse} that
\begin{align*}
&|\int \theta_\varepsilon (d\rho_t-d\rho_0)|
\le |t| \sup_\tau | \bigl ((1-\LL_\tau)^{-1}[\divv (X_\tau)\rho_\tau + \langle  \gradd \rho_\tau, X_\tau \rangle ]\bigr )  (\theta_\varepsilon)|\\
&\qquad\qquad \le \widetilde C_X C' |t|  \sup_\tau \|\rho_\tau\|_{W^{u,s}_{p'}} \|\theta_\varepsilon\|_{H^{1+\tilde r}_p}
\le C \cdot C' \cdot \widetilde C_X \cdot  |t|  \sup_\tau \|\rho_\tau\|_{W^{u,s}_{p'}} \|\theta\|_{H^{r}_p}  \varepsilon^{r-1-\tilde r}\, .
\end{align*}
(We used $\|\theta_\varepsilon\|_{(W^{u-1,s-1}_{p'})^*}\le C' \|\theta_\varepsilon\|_{H^{1+\tilde r}_p}$
for any $\tilde r >|s|$.)

On the other hand,
\begin{equation}\label{easy0}
|\int (\theta_\varepsilon-\theta) d\rho_t|+
|\int (\theta_\varepsilon-\theta) d\rho_0|
\le 2 \|\theta_\varepsilon - \theta\|_{H^{\tilde r}_p}\max_{\tau=0,t} \|\rho_\tau\|_{H^{-\tilde r }_{p'}}\, .
\end{equation}
Since $\|\rho_\tau\|_{H^{-\tilde r}_{p'}}\le \|\rho_\tau\|_{W^{u,s}_{p'}}$
if $\tilde r>|s|$, we have $\max_{\tau=0,t} \|\rho_\tau\|_{H^{-\tilde r}_{p'}}<\infty$.
Since $\tilde r<r<1/p$, we have 
\begin{equation}\label{easy} \|\theta_\varepsilon - \theta\|_{H^{\tilde r}_p} \le C \varepsilon^{r-\tilde r}\|\theta\|_{H^{r}_p}\, ,
\end{equation}
 e.g., by \cite[Lemma 5.4]{BLi}. 
Choose $\varepsilon=|t|$.  For each $0<\alpha<r$, we can get
$$
|t|^{1+r-1-\tilde r}=|t|^{r-\tilde r}\le |t|^\alpha \, , 
$$
by taking small enough   $\tilde r>|s|>0$. This proves   that
$t \mapsto \int \theta \rho_t$ is $\alpha$-H\"older.
\end{proof}

\section{Linear response for hyperbolic attractors and observables with transversal singularities}
\label{sec2}

As pointed out in Remark~\ref{basic}, a necessary condition for our study is
that $\int \theta \, d\rho_t$ be well defined for the observable $\theta$. In order
to get linear response, the formula \eqref{LRF} in the proof of Theorem~\ref{main}
shows that we  need  the stronger condition that $\theta$ belongs to the dual
of a space containing  the distribution
$\langle \gradd \rho_t, Y\rangle$ for a smooth
vector field $Y$, in particular containing the components (in the sense of distributions) of $\gradd \rho_t$. We may thus hope
that observables $\theta$ with singularities satisfy this stronger condition if
their ``wave front set'' (WFS)
is transversal to the WFS of $\gradd \rho_t$. (We only use the  notion of a WFS in this informal discussion, and the
reader will  not need to know its precise definition, which is not exactly the classical one since
our objects  are $C^r$ for  some  finite $r >1$, possibly small with respect
to the dimension. It suffices to
mention that, in the present setting, the WFS of $\rho_t$ and $\gradd \rho_t$ are contained
in stable cones.)

Let us be more concrete, considering  the anisotropic spaces $W^{u,s}_p$ (containing $\rho_t$) used in the proof of
Theorem~ \ref{main}. To fix ideas, we work in dimension two and focus on the toy-model situation
of charts $(x,y)$
where the horizontal direction $y=0$ is expanding and the vertical direction $x=0$ is contracting.
We consider $p=p'=2$ to simplify matters, and in this informal discussion
 we shall pretend that the dual of $W^{u,s}_2$ is $W^{-u,-s}_{2}$. Consider 
the toy function $\theta(x,y)=\Theta(x) h(y)$, where $\Theta$ is the Heaviside function, and $h$ is a
$C^\infty$ function with rapid decay (we ignore compact support issues in the $x$ variable
in the present heuristic  outline). This is clearly the simplest  transversal 
toy function. In Fourier space, with parameters $(\xi,\eta)$, we find 
$\hat \Theta(\xi) \hat h(\eta)$, where
the decay
 of $\hat \Theta$ at infinity is  $\sim |\xi|^{-1}$ and for any
$N\ge 1$ the decay of $\hat h(\eta)$ at infinity is  $\sim C_N |\eta|^{-N}$. 
 It is thus
 reasonable to hope that a function of the
type of $\theta(x,y)$ belongs not only  to the dual $W^{-u,-s}_{2}$ of a space
 $W^{u,s}_2$  for some
$s<0<u$, but also 
to the dual $W^{-(u-1),-(s-1)}_{2}$ of a space $W^{u-1,s-1}_2$.
Indeed , in  
 any fixed ``stable'' cone $|\xi|\le \tilde c |\eta|$,  we can get the required decay  $C (1+|\xi|^2+|\eta|^2)^{(s-1)/2}$ (with $s<0$ so that $s-1<-1$) by choosing $N$ large enough depending
on $\tilde c$. The condition in the ``unstable'' cone is satisfied if $(u-1+1)2>2$, that is, if
$u>1$.
(We would like to emphasize that, although we have just argued that
functions of the type $\theta(x,y)$ may
lie in anisotropic spaces $W^{u,s}_p$, they are in general not bounded multipliers in
these spaces, see \cite{baladijoel}.)

Our second main theorem below (Theorem~\ref{main2}) will justify the optimism induced by the above toy-model
computation. 

To formulate this second result, we first recall the notion of admissible stable leaves from 
Gou\"ezel and Liverani \cite[\S 3]{GL1}.
 Let $t\mapsto f_t$, for $t\in [-\epsilon_0,\epsilon_0]$, be a $C^3$ family
of  $C^4$  diffeomorphisms $f_t$   on a smooth $d$-dimensional Riemann manifold $M$, with a transitive compact 
hyperbolic
attractor $\Lambda_t\subset M$ (as defined in the beginning
of Section~\ref{sec1}), and let $V$ be a common attracting neighbourhood for all
the $f_t$ (up to reducing $\epsilon_0$). We may and do assume that $\overline V$ is compact.
Using an adapted Mather metric and further reducing $\epsilon_0$ if necessary, we can assume that there exists $0<\nu<1$ so that
for all small enough $|t|$, the expansion of $D_xf_t$ along $E^u_x$
 is stronger than $\nu^{-1}>1$, while its contraction
along $E^s_x$ is stronger than $\nu<1$, and the angle
between $E^s_x$ and $E^u_x$  is everywhere arbitrarily
close to $\pi/2$. For small enough $\kappa>0$, we define the stable
cone at $x\in V$ by
  \begin{equation*}
  \CC^s(x)=\left\{ w+v \in T_x M \mid w\in E^s(x)\, , v \perp E^s(x)\, , \|v\|
  \leq \kappa \|w\| \right\} \, .
  \end{equation*}
If $\kappa>0$ is small enough then
$D_xf_0^{-1}(\CC^s(x)\backslash \{0\})$ belongs to
the interior of $\CC^s(f_0^{-1}(x))$, and $D_xf_0^{-1}$ expands the
vectors in $\CC^s(x)$ by $\nu^{-1}$.

\begin{definition}[Admissible charts]\label{notat0}
There exist 
real numbers $r_i \in (0,1)$ and $C^4$ coordinate charts
$\psi_1,\ldots,\psi_N$, with $\psi_i$  defined on 
$(-r_i,r_i)^d\subset \real^d$ (with its Euclidean norm), such
that the attracting neighbourhood $V$ is covered by the open sets
$$\bigl(\psi_i((-r_i/2,r_i/2)^d)\bigr)_{i=1\ldots N}\, , $$
and the following conditions hold: $D\psi_i(0)$ is an isometry,  $D\psi_i(0)\cdot \bigl(\real^{d_s}\times\{0\}\bigr)= E^s(\psi_i(0))$,
 and the $C^4$-norms of $\psi_i$ and its inverse are bounded by
$1+\kappa$ (with $\kappa$ as introduced above).
\end{definition}

 Next, we may pick $c_i\in (\kappa,2\kappa)$
such that the corresponding stable cone in charts 
$$\CC^s_i=\{ w+v\in \real^d
\mid w\in \real^{d_s}\times \{0\}, v\in \{0\}\times \real^{d_u}, \|v\|\leq
c_i \|w\|\}
$$  
satisfies   $D_x\psi_i (\CC^s_i) \supset \CC^s(\psi_i (x))$ and $D_{\psi_i(x)} f_t^{-1}(
D_x\psi_i(x) \CC^s_i) \subset \CC^s(f_t^{-1}( \psi_i(x)))$
for any  $x\in
(-r_i,r_i)^d$ and all $|t|<\epsilon_0$ (reducing $\epsilon_0$ if necessary).

Let $G_i(K)$ be the set of graphs of $C^4$ maps $\chi: U_\chi \to (-r_i,r_i)^{d_u}$ defined on a
subset $U_\chi$ of $(-r_i,r_i)^{d_s}$, with 
$$|D\chi| <  c_i \mbox{ and } |\chi|_{C^4} \le K \, ,$$ 
where a
function defined on a subset $U$ of $\real^{d_s}$ is 
$C^4$ if it admits a $C^4$ extension to an open
neighbourhood of $U$,  the norm being the infimum of the norms of such
extensions. 
In particular, the tangent space to the graph of $\chi$ belongs to
the interior of the cone $\CC^s_i$.

Uniform hyperbolicity of $f_0$ implies (see \cite[Lemma 3.1]{GL1}) that
 if $K$ is large enough, then there exists $K'<K$
such that, for any $W\in G_i(K)$ and any $1\le j\le N$, the
set $\psi_j^{-1}(f_t^{-1} (\psi_i(W)))$ belongs to $G_j(K')$
for all $|t|<\epsilon_0$ (reducing again $\epsilon_0$ if needed).

\begin{definition}[Admissible graphs and admissible stable
    leaves]\label{notat}
    An admissible graph is  a map $\chi$ defined on a
ball 
$$\ovB(w,K_1\delta)\subset (-2r_i/3,2r_i/3)^{d_s}
$$ 
for some small enough $\delta>0$
and
large\footnote{If $\kappa>0$ chosen above is small enough, then
$\nu (1+\kappa)^2 \sqrt{1+4\kappa^2}<1$. Let $K_1>1$ be such
that
  $K_1>1+K_1 \nu_1 (1+\kappa)^2\sqrt{1+4\kappa^2}$, and
take $\delta>0$  so that $6 K_1\delta < \min_i(r_i)$, see \cite{GL1}.} 
enough $K_1$, taking
its values in $(-2r_i/3,2r_i/3)^{d_u}$
 with
$\operatorname{range}(\Id,\chi)\in G_i(K)$. 
The set of admissible stable leaves
is 
\begin{align*}
 \bigl\{ W:=\psi_i \circ (\Id,\chi)(\ovB(w,\delta)) \mid&
  \chi:\ovB(w,K_1\delta) \to \real^{d_u} \\
  &\qquad\text{ is an admissible graph on } B_i:=(-2r_i/3,2r_i/3)^{d_s}
  \bigr\}\, .
\end{align*}
\end{definition}

We  may now state our second theorem, which
exhibits   a simple  transversality condition  which is sufficient to ensure linear response
for hyperbolic attractors: 

\begin{theorem}
[Linear response for transversal Heaviside singularities]
\label{main2} Let $t\mapsto f_t$, for $t\in [-\epsilon_0,\epsilon_0]$, be a $C^3$ family
of  $C^4$ 
diffeomorphisms $f_t$   on a smooth Riemann manifold $M$, with a transitive compact 
hyperbolic
attractor $\Lambda_t\subset M$. Let $\rho_t$ be
the (unique) SRB measure of $f_t$ on $\Lambda_t$.  
Let
$$
\theta(x)=h(x) \Theta(g(x)-a)\, ,
$$
where $\Theta:\real\to \real$ is the Heaviside function, 
$h:M \to \complex$ is $C^3$ and  $g:M \to \real$ is $C^4$, with
$a\in \real$ not a critical value of $g$, and so
that\footnote{This assumption implies the  announced transversality condition: The normal to $\WW_a$
at each point lies in the unstable cone.}
 the intersection 
 $$\WW_a:=\{ x\mid g(x)=a\} \cap \supp(h)$$ 
admits
  a $C^4$ foliation by admissible stable leaves. Then there exists $\epsilon_1\in (0,\epsilon_0]$ so that for 
any $|t_0|<\epsilon_1$ the map
$
t\mapsto \int_M \theta\,  d\rho_t
$
is  differentiable at $t_0$, with (recalling \eqref{defLt} and \eqref{X_t})
\begin{align}
\partial_t (\int_M \theta \, d \rho_t )|_{t=t_0}
\label{LRF2} & =-
\bigl ( (1-\LL_{t_0})^{-1}[\divv (X_{t_0})\rho_{t_0}+ \langle \gradd \rho_{t_0}, X_{t_0}\rangle ]\bigr ) (\theta)
\, , 
\end{align}
where the gradient is in the sense of distributions.

In addition, if $(f_{t_0}, \rho_{t_0})$
is\footnote{As usual, we can reduce to the mixing case by spectral decomposition (see e.g. \cite[Theorem 18.3.1]{KH}).} mixing, we have the exponentially convergent sums
\begin{align}
\label{FDT}
\partial_t (\int_M \theta \, d \rho_t )|_{t=t_0}
&=-\sum_{k=0}^\infty
\bigl (  \LL^k_{t_0} [\divv (X_{t_0})\rho_{t_0}+ \langle \gradd \rho_{t_0}, X_{t_0}\rangle ]\bigr ) (\theta)\\
\label{LRF3}&=\sum_{k=0}^\infty \rho_{t_0}\bigl ( \langle \gradd (\theta \circ f_{t_0}^k), X_{t_0}\rangle 
\bigr ) \, .
\end{align}
\end{theorem}

If the stable dimension $d_s$ is equal to $d-1$ (for example, if
$d=2$ with $d_s=d_u=1$)  
then the assumption 
reduces to requiring that 
 $\WW_a=\{ x\mid g(x)=a\} \cap \supp(h)$ 
 is an  admissible stable leaf.

More general transversal test functions can be considered: Our proof of Theorem~\ref{main2} actually
shows that any function $\theta$ lying in the dual  to
a Gou\"ezel--Liverani space $\BB^{0,|s-2|}$ (the corresponding norm is recalled in \eqref{st1})
with $s<0$ satisfies the conclusion of the theorem. Since our goal here is to highlight
a sufficient condition which works for a natural class of test functions, we refrain from
a more general statement in order to keep the proof as free as possible from technicalities.
(In any case, the regularity conditions in the blueprint theorem above are probably not optimal, in
particular we expect that lying in the dual  to
a  space $\UU^{0,|s-1|}_p$ with $s<0$ from \cite{baladijoel} should suffice.)

\begin{remark}[Potential generalisations involving the Dirac mass]\label{blueprint}We expect that there
is an analogue of Theorem~\ref{main2} when
$h(x)\Theta(g(x)-a)$ is replaced by the distribution $h(x)\delta_D(g(x)-a)$ 
if  $a$ is not critical and  $\WW_a$  is foliated by
admissible stable leaves, up to increasing the smoothness requirement.
(See also \cite[Prop 4.4]{GL1} for an occurrence of the Dirac mass in the anisotropic setting.)
\end{remark}

\begin{proof}[Proof of Theorem~\ref{main2}]
 To prove Theorem~\ref{main2}, it is
convenient  to use the anisotropic Banach spaces from
\cite{GL1, GL2}, based on taking suprema over admissible stable leaves (instead
of working with cones in Fourier space as we did with the spaces from \cite{BTs}
used in the proof of Theorem~\ref{main}).

Since each $f_t$ is $C^4$, the transfer operator $\LL_t$
has a spectral gap  on the anisotropic spaces $\BB^{u,|s|}$ of \cite{GL1} for integer $u\ge 1$
and real $s<0$ with $u+|s|<3=4-1$.  Using
 admissible stable leaves and the notations
introduced from Definition~ \ref{notat}, we  recall the definition of
the norm of $\BB^{u,|s|}$ in coordinates (see \cite[Lemma 3.2]{GL1}). Let 
$u\in \{ 1, 2\}$, let
 $-3+u<s<0$  be a real number, and set
(the test function $\omega$ below is compactly supported in $\ovB(w,\delta)$, 
the measure $dm$ is  Lebesgue
measure on $\real^{d_s}$,
and $\chi$ ranges over admissible graphs on $B_i$)
  \begin{equation}\label{st1}
  \|\varphi\|_{u,s} =  \max_{\substack{0 \le u'\le u\\ u'\in \integer}}\,  \max_{\substack{|\alpha|=u'\\1\le i\le N}}\;
  \sup_{\chi  }\,\, 
  \sup_{ |\omega|_{C^{|s|+u'}} \le 1}\,\,  \int_{B(w,\delta)} 
  [\partial^\alpha( \varphi \circ \psi_i)] \circ (\Id,\chi) \cdot \omega \, dm\, .
  \end{equation}
  The space $\BB^{u,|s|}$ is then the closure of $C^3$ functions for the above norm.

 We need further preparations:
We may
and shall focus on a neighbourhood
of $\WW_a$, since the function $\theta$ is $C^2$ elsewhere.  
Since $a$ is not a critical value, $\WW_a$ is a codimension-one manifold
by the implicit function theorem. Our assumption guarantees that this manifold
has a $C^4$ foliation by $d_s$-dimensional leaves which are stable admissible
leaves.
The idea is to foliate a chart neighbourhood by stable leaves, e.g. like in the proof
of \cite[Lemma 2.11]{BDL}. More precisely, we consider the $C^4$ foliation of $\WW_a$ by stable
leaves (which exists by assumption), noting $W_a(z)$ the leaf through $z\in \WW_a$. Fixing  small balls around finitely many points $x_j\in \WW_a$, we
may use $C^4$ charts $\tilde \psi_j$ satisfying the assumptions  in
Definition~\ref{notat0} and, in addition, mapping    $\real^{d_s}\times \{\vec 0\}$  to 
$W_a({x_j})$, while mapping the codimension-one hyperspace $\real^{d-1}\times \{0\}$  to $\WW_a$. 
(Recall also that the $d$-th axis $\{0\} \times \real$ is transversal to the stable cone $\CC^s_j$.) 
It will be convenient to require the charts to satisfy the following stronger assumption:
First notice that since $a$ is not a critical value of $g$, there is a small neighbourhood $I_a$ of
$a$ so that any $b\in I_a$ is not a critical value of $g$. In addition,
$\WW_b$ is also foliated by stable leaves. In fact, we may ensure that there exist
$C^3$ local diffeomorphisms $\tau_j:\real\to \real$ with $\tau_j(0)=0$,
$\tau_j'(0)=1$, and $\tau_j''(0)=0$, so that for any
$b \in I_a$, the chart
 $\tilde \psi_j$  maps    $\real^{d_s}\times \{\tau_j(a-b)\}$  to 
$W_b({x_j})$ and  maps the codimension-one hyperspace $\real^{d-1}\times \{\tau_j(a-b)\}$  to $\WW_b$. 
We
let $\tilde \phi_j$ be an adapted partition of unity.

\medskip
We now move on to the proof. Assume without loss of generality that $t_0=0$.
Set
\begin{equation}\label{gma}
\gamma:=(1-\LL_{0})^{-1}[\divv (X_{0})\rho_{0}+ \langle \gradd \rho_{0}, X_{0}\rangle ]\, .
\end{equation}
Since the results of \cite{GL1} imply that $\rho_t\in \BB^{u, |s|}$ for $u=2$ and $-1<s<0$,
and since $X_0$ is a $C^2$ vector field,
  we have $\gamma\in \BB^{u-1, |s-1|}\subset \BB^{u-2, |s-2|}$.
 (The Leibniz inequality is not difficult here since each $u'$ in \eqref{st1} is an integer, in particular,
we do not need different sets of ``cones'' as in the proof of Theorem~\ref{main}.)

We can apply \cite[Theorem 2.7, \S 8]{GL1}  (as we are using the same
spaces here): The result there says that
the map $t \mapsto \rho_t \in \BB^{u-2, |s-2|}$ 
is differentiable and, for any $\varphi$ in the dual space to $\BB^{u-2, |s-2|}$
we have (recalling the notation \eqref{gma})
\begin{equation}\label{limm}
\partial_t \int \varphi \, d \rho_t|_{t=0}=-\gamma(\varphi)=-
\bigl( (1-\LL_{0})^{-1}[\divv (X_{0})\rho_{0}+ \langle \gradd \rho_{0}, X_{0}\rangle ]\bigr ) 
(\varphi)\, .
\end{equation}
(Contrary to the situation in the proof of Theorem~\ref{main}, it is not clear that
we can take
$\varphi$ in the dual space to $\BB^{u-2, |s'-1|}$ for $s-1<s'<s<0$ in \eqref{limm}.)

To prove the first claim of the theorem, it thus suffices to show that  the function
$\theta(x)=h(x) \Theta(g(x)-a)$ belongs to  the dual space to $\BB^{u-2, |s-2|}$. 
(Recall that $u=2$ and $-1<s<0$.)
The special properties of the charts $\tilde \psi_j$ derived from
our transversality assumption will be essential here.
Let $\tilde \gamma$ be a $C^3$ function on $M$ supported in a small ball around $x_j$
and   put
$$h_j=|\det D\tilde \psi_j|(\tilde \phi_j \cdot h)\circ \tilde \psi_j\, .
$$ 
Then the contribution of
$\int \tilde \gamma \theta \, dm=\int \tilde \gamma(y) h(y)\Theta(g(y)-a) \, dm$ 
to the $j$th chart is
\begin{align*}
&\int_{\real^d} |\det D\tilde \psi_j|(\tilde  \phi_j \theta \tilde \gamma) \circ \tilde \psi_j \, dx\\
&\qquad\qquad\qquad\qquad=
\int_{[-1,1]^{d_u}} \Theta(x_d)   dx_{d_s+1}\ldots dx_{d} \int_{B_j} dx_1\ldots dx_{d_s} h_j(x) \tilde \gamma (\tilde \psi_j(x))\\
&\qquad\qquad\qquad\qquad=
\int_{[0,1]}  dx_d \int_{[-1,1]^{d_u-1}} dx_{d_s+1}\ldots dx_{d-1} \int_{B_j} dx_1\ldots dx_{d_s} h_j(x) \tilde \gamma (\tilde \psi_j(x)) \, .
\end{align*}
Obviously,
\begin{align*}
&| \int_{[0,1]}  dx_d \int_{[-1,1]^{d_u-1}} dx_{d_s+1}\ldots dx_{d-1} \int_{B_j} dx_1\ldots dx_{d_s} h_j(x) \tilde \gamma (\tilde \psi_j(x))| \\
&\qquad\qquad\qquad\qquad\qquad\qquad\le
2^{d_u-1} \sup_{x_{d_s+1}, \ldots, x_d}
|\int_{B_j} dx_1\ldots dx_{d_s} h_j(x) \tilde \gamma (\tilde \psi_j(x))|
\, .
\end{align*}
Since $h_j$ is $C^3$, while $s-2 >-3$ and $u =2$, taking $u'=0$ in the
definition \eqref{st1} of the norm, the integral
over $B_j\subset [-1,1]^{d_s}$ in the right-hand side is bounded by  $\|h_j\|_{C^2} \|\tilde \gamma\|_{u-2,s-2}$.  In particular,\footnote{In practice, we can take the sequence corresponding
to $(1-\LL_0)^{-1} [\divv X_0 \LL_0^n(1_V)+\langle X_0, \gradd (\LL_0^n(1_V)\rangle ]$.}
we can let $\tilde \gamma$ tend to  $\tilde \phi_j \cdot \gamma \in \BB^{u-2, |s-2|}$. We have proved  that $\theta$  lies in
the dual space to $\BB^{u-2,|s-2|}$.

It remains to show  \eqref{FDT}--\eqref{LRF3} in the mixing case.
For this, first notice that (see \cite{GL2})
mixing implies that $1$ is the only eigenvalue of $\LL_{0}$
on the unit circle so that, using \eqref{boundaryless}, 
the distribution $(1-\LL_{0})^{-1}(\gamma)$ may be replaced
by the exponentially convergent --- in the Banach norm ---
sum $\sum_{k=0}^\infty \LL_{0}^k(\gamma)$. Second, note that
$(\LL_{0}^k(\gamma))(\theta)=\gamma(\theta\circ f_{0}^k)$.
Finally, integrate by parts, in the 
sense of distributions,
noting
that the derivative of $\theta$ is  not only a distribution of order zero
but also an element of the dual of $\BB^{u,|s|}$. This can be checked by using the
ideas developed above, and is left to the reader.
\end{proof}

\appendix
\section{Definition of the anisotropic spaces $W^{u,s}_p$ from \cite{BTs}.}
\label{appA}

For the reader's convenience, we recall the
definition of the space $W^{u,s}_p(T,V)=W^{u,s}_{p,\dagger}(T,V)$ from \cite[App. A]{BTs}
associated to a $C^\kappa$ diffeomorphism $T$ with a transitive compact hyperbolic attractor $\Lambda$
with an attracting neighbourhood $V$, for $\kappa>1$.
(Formally the paper \cite{BTs} only considers the Anosov case when $\Lambda=M$, but, since $E^u$ and $E^s$ can
be extended to  $\overline V$, in the sense of \cite{BT2}, this does not cause
any problems.)

A cone  in $\real^d$ is a subset 
which is invariant under scalar
multiplication.
We write 
$\cone \Subset \cone'$  for two cones $\cone$ and $\cone'$
in $\real^d$
if
$\overline \cone\subset \mbox{ interior} \, (\cone' )\cup \{0\}$.
A cone $\cone$ is called $d'$-dimensional 
 if $d'\ge 1$
is the maximal dimension of a linear subset of $\cone$.

\begin{definition}[Cone pairs and cone systems] A cone pair is 
$\cone_\pm=(\cone_+,\cone_-)$, 
where  $\cone_+$ and $\cone_-$ are
closed cones in $\real^d$,  with nonempty interiors, of respective
dimensions $d_u$ and $d_s$ and so that $\cone_+\cap\cone_-=\{0\}$.
A cone system is 
$$\Theta=(\cone_\pm,\varphi_+,\varphi_-)\, ,  $$  
with $\cone_\pm=(\cone_+,\cone_-)$ a cone pair and
$\varphi_\pm:\sphere\to [0,1]$ two
$C^{\infty}$ functions on the unit sphere $\sphere$ in $\real^d$ satisfying
$$
\varphi_+(\xi)=
1 \mbox{ if $\xi\in \sphere\cap \cone_{+}$,}\quad
\varphi_+(\xi)=0 \mbox{ if $\xi\in \sphere\cap \cone_{-}$,} \quad \varphi_-=1-\varphi_+\, .
$$
\end{definition}

\medskip

Next, introduce for real numbers $u$ and $s$  the functions 
\[
\Psi_{u,\Theta_+}(\xi)=(1+\|\xi\|^2)^{u/2}\varphi_+
\left (\frac{\xi}{\|\xi\|}\right )\quad\mbox{and}
\quad \Psi_{s,\Theta_-}(\xi)=(1+\|\xi\|^2)^{s/2}
\varphi_-\biggl (\frac{\xi}{\|\xi\|} \biggr )\, .
\]
For a cone system $\Theta$
and $\varphi\in C^\infty_0(\real^d)$,
we define norms  
for $1<p<\infty$ and $ s \le 0 \le u$ by
\begin{eqnarray}
\label{daggernorms}  \|\varphi\|_{W^{\Theta, u,s}_{p}}
&=&\|\Psi_{u,\Theta_+}^{Op}(\varphi)\|_{L_p}+
\|\Psi_{s,\Theta_-}^{Op}(\varphi)\|_{L_p}\, ,
\end{eqnarray}
where for a  function $\Psi: \real^d \to \real$ the operator $\Psi^{Op}$ is defined by the convolution
\begin{equation*}
\Psi^{Op} (\varphi)=(\FFF^{-1} \Psi) * \varphi\, , 
\end{equation*}
 where  $\FFF^{-1}$  denotes the inverse Fourier transform
\begin{align*}
\FFF^{-1}(\Psi)(x)&= \frac{1}{(2\pi)^d}
\int_{\real^d} e ^{i x\xi}  \Psi(\xi) \, d \xi \, , \quad x\in \real ^d \, .
\end{align*}

\begin{definition}[Admissible charts and partition of unity]
Admissible charts and partition of unity
for $T$ are:
A finite system of $C^{\infty}$ local charts $\{(V_i, \psi_i)\}_{1\le N}$, with 
open subsets
$V_i\subset M$, and diffeomorphisms
$\psi_i : U_ i\to V_i$ such that $\overline V \subset \cup_i V_i$, and  
$U_i$  
is a bounded open subset of $\real^d$ for each $1 \le i\le N$, together
with a  finite $C^{\infty}$ partition of  unity  $\{\phi_i\}$ for $M$,
   subordinate to the cover $\VV=\{V_i\}$ of $\overline V$. 
\end{definition}

\begin{definition}[Admissible cone pairs]
Since $\Lambda$ is a hyperbolic attractor for $T$ we may choose  $C^{\infty}$ 
admissible local charts $\{(V_i, \psi_i)\}_{1 \le i\le N}$
and cone pairs
$
\{\cone_{i,\pm}=(\cone_{i,+}, \cone_{i,-})\}_{1 \le i\le N}
$, 
so   that the following conditions hold\footnote{$\cone_{i,\pm}$ are  locally constant cone fields
in the cotangent bundle $T^* \real^d$, so that the conditions are expressed
with respect to normal subspaces.}:
\begin{itemize}
\item If $x\in V_i$, the cone
$(D\psi_i^{-1})^{*}_x(\cone_{i,+})$    
contains the ($d_u$-dimensional) normal subspace of $E^s(x)$, and the cone
$(D\psi_i^{-1})^{*}_x(\cone_{i,-})$
contains the ($d_s$-dimensional)
normal subspace of  $E^u(x)$.
\item 
If  $T(V_i)\cap V_{i'}\ne \emptyset$, then the $C^\kappa$ map 
\[
F=\psi^{-1}_i
\circ T^{-1}\circ \psi_{i'}
:\psi_{i'}^{-1}(T(V_i)\cap V_{i'}) \to U_i \, ,
\]
extends to a  bilipschitz $C^1$
diffeomorphism of $\real^d$ 
so that, using $A^{tr}$ to denote the transposition of a matrix $A$,
$$DF_{x}^{tr}(\real^d\setminus \cone_{i,+}) \Subset \cone_{i',-}\, , \qquad 
\forall x\in \real^d\, .
$$
\end{itemize}
\end{definition}

\bigskip

\begin{definition}[The spaces $W_{p}^{u,s}(T,V)$]
Fix a finite set of admissible charts, a partition of unity, and
admissible cone pairs, and choose a cone system for each cone pair. Then for any
$1<p<\infty$ and $u-(\kappa-1)< s \le 0 \le u$,
the Banach space $W_{p}^{u,s}(T,V)$ is
the completion of $C^{\infty}(\overline V)$ for the  norm
\begin{equation}
\|\varphi\|_{W_{p}^{u,s}}:=\max_{1 \le i \le N} 
\|(\phi_i\cdot \varphi)\circ \psi_i^{-1}\|_{W^{\Theta_i, u,s}_{p}}\, ,
\end{equation} 
recalling \eqref{daggernorms}
(We do not claim that choosing different charts, partitions of unity, cone pairs, or cone systems
would lead to equivalent norms.) 
\end{definition}

\end{document}


\newcommand{\comments}[1]{\marginpar{\footnotesize #1}} 

\newtheorem{proposition}{Proposition}[section]
\newtheorem{lemma}[proposition]{Lemma}
\newtheorem{sublemma}[proposition]{Sublemma}
\newtheorem{theorem}[proposition]{Theorem}

\newtheorem{maintheorem}{Main Theorem}
\newtheorem{corollary}[proposition]{Corollary}

\newtheorem{ex}[proposition]{Example}

\theoremstyle{remark}

\newtheorem{remark}[proposition]{Remark}

\theoremstyle{definition}
\newtheorem{definition}[proposition]{Definition}
\newcommand{\ovB}{\bar{B}}
\def\Erg{\mathrm{Erg}\, }
\def\cone{\mathbf{C}} 
\def\real{\mathbb{R}}
\def\sphere{\mathbf{S}^{d-1}}
\def\integer{\mathbb{Z}}
\def\complex{\mathbb{C}}
\def\BBB{\mathbb{B}}
\def\supp{\mathrm{supp}}
\def\var{\mathrm{var}}
\def\sgn{\mathrm{sgn}}
\def\sp{\mathrm{sp}}
\def\id{\mathrm{id}}
\def\Imm{\mathrm{Image}}
\def\cc{\Subset}
\def\D{\mathrm {d}}
\def\I{i}
\def\E{e}
\def\Lip{\mathrm{Lip}}
\def\BB{\mathcal{B}}
\def\CC{\mathcal{C}}
\def\DD{\mathcal{D}}
\def\EE{\mathcal{E}}
\def\FF{\mathcal{F}}
\def\GG{\mathcal{G}}
\def\II{\mathcal{I}}
\def\JJ{\mathcal{J}}
\def\KK{\mathcal{K}}
\def\LL{\mathcal{L}}
\def\LLL{\mathbb{L}}
\def\MM{\mathcal{M}}
\def\NN{\mathcal{N}}
\def \OO{\mathcal {O}}
\def \PP{\mathcal {P}}
\def \QQ{\mathcal {Q}}
\def \RR{\mathcal {R}}
\def\SS{\mathcal{S}}
\def\TT{\mathcal{T}}
\def\UU{\mathcal{U}}
\def\VV{\mathcal{V}}
\def\YY{\mathcal{Y}}
\def\ZZ{\mathcal{Z}}
\def\FFF{\mathbb{F}}
\def\PPP{\mathbb{P}}

\title[Corrigendum to linear and fractional response]{Corrigendum to: Linear and fractional response for the SRB measure of smooth hyperbolic attractors and discontinuous observables}

\author{Viviane Baladi, Tobias Kuna, and Valerio Lucarini} 

\date{\today}

\begin{abstract}
The first main result of \cite{BKL} is modified as follows:
 For any  $\theta$ in the  Sobolev space
$H^r_p(M)$, with $1<p<\infty$ and  $0<r<1/p$, the map $t\mapsto \int \theta\, d\rho_t$ is $\alpha$-H\"older continuous for
all  {\boldmath$\alpha <r- |\log \JJ|/(p|\log \nu_s|)$}
 {\bf where
 {\boldmath$\JJ\le 1$} is the strongest volume contraction and  {\boldmath$\nu_s<1$} is the weakest contraction.} This  applies to
 $\theta(x)=h(x)\Theta(g(x)-a)$ (for all 
 {\boldmath$\alpha <1- |\log \JJ|/|\log \nu_s|$}) for $h$ and $g$
 smooth and $\Theta$ the
 Heaviside function, if $a$ is not a critical value of $g$.''
 \end{abstract}
\thanks{We thank S. Galatolo for pointing out the counter-example to the initial statement
of Theorem ~2.1. We are grateful to S. Gou\"ezel, P.-A. Guih\'eneuf,  and Beno\^\i t Saussol
for useful comments.}


\maketitle

\section{Correcting \cite[Theorem 2.1]{BKL} on fractional response}

The essential spectral radius of the operator $\LL_t$ defined by \cite[(2.1)]{BKL} acting
on $W^{u,s}_{p',\dagger}(f_t,V)$ for $s<0<u$ is only guaranteed to be strictly smaller than $1$ by \cite{BT1} if  
\begin{equation} \label{star_0} \lim_{m\to \infty} 
\sup_{\Lambda_t} [|\det D_xf^m_t|^{-1/p}
\max\{(\nu_s(f^m_t))^{|s|},( \nu_u(f^m_t))^{-u}\}]^{1/m}<1\, ,
\end{equation}
where 
$\nu_s(x,f_t^{m})=\sup_{v\in E^s(x)\setminus\{0\}}
\frac{\|D_xf_t^{m}(v)\|}{\|v\|}$,
 $\nu_u(x,f_t^{m})
=\inf_{v\in E^u(x)\setminus\{0\}}
\frac{\|D_xf_t^m(v)\|}{\|v\|}$  for $x\in \Lambda_t$ and $m\in \integer_+$.
Because of this, \cite[Theorem 2.1]{BKL} must  be modified:

\begin{theorem}[Fractional response]
Fix $\beta \in (0,1)$. Let $t\mapsto f_t$, for $t\in [-\epsilon_0,\epsilon_0]$, be a $C^{2+\beta}$ family
of  $C^3$  diffeomorphisms $f_t$   on a smooth Riemann manifold $M$, 
so that $f_t$ has a transitive compact 
hyperbolic
attractor $\Lambda_t\subset M$. Let $\rho_t$ be
the (unique) SRB measure of $f_t$ on $\Lambda_t$.  
{\bf Assume that, for some {\boldmath $1<p<\infty$}  and all} {\boldmath $t$},\footnote{{\bf If each {\boldmath $f_t$ is $C^{N+1}$} and the family is $C^{N+\beta}$
for \boldmath{$N\ge 3$}, then \eqref{star_1}  can be replaced by the condition
{\boldmath $\nu_{u,t}^{-(N+\beta-2)}<\JJ_t^{\frac{1}{p}}$}.}}
{\boldmath \begin{equation}
\label{star_1}
\nu_{u,t} ^{-\beta}:=\lim_{m\to \infty} 
\sup_{\Lambda_t} [ \nu_u(x,f^m_t)]^{-\frac{\beta}{m}}<\JJ_t^{\frac{1}{p}}:= \lim_{m\to \infty} \inf_{\Lambda_t} |\det Df^m_t(x)|^{\frac{1}{pm}}\, .
\end{equation}}
Let $\theta :M \to \complex$,  be so that\footnote{{\boldmath$H^r_p(\bar V)$} {\bf is the closure of {\boldmath$C^\infty(\bar V)$} for the norm of} 
$H^{r}_p(M)$.}
 {\boldmath$\theta \in H^r_p(\bar V)$} for  some
 $0< r <1/p$. Then there exists $\epsilon_1\in (0, \epsilon_0]$ so that  
$t\mapsto \int_M \theta\, d\rho_t$
is $\alpha$-H\"older continuous on $[-\epsilon_1, \epsilon_1]$
for any
{\boldmath \begin{equation}
\label{star_2}  0 <\alpha <r-\sup_{|t|<\epsilon_1} 
\frac{1}{p}\frac{|\log \JJ_t|}{|\log \nu_{s,t}|}\, .\end{equation}} 
\end{theorem}

{\bf If 
 {\boldmath $r$} is close enough to {\boldmath $1/p$}, the r.h.s. of 
\eqref{star_2} is $>0$ if, for\footnote{{\bf Use semi-continuity.}} all $t$}
{\boldmath \begin{equation}\label{star_3} 
\nu_{s,t}:=\lim_{m\to \infty} \sup_{\Lambda_t} [\nu_s(x,f^m_t)]^{\frac{1}{m}}
<\JJ_t\, .
\end{equation}}

Remark 2.3 from \cite{BKL} should be modified correspondingly: If\footnote{{\bf By Footnote 1,
condition \eqref{star_1} can be waived  if   \boldmath{$g, h, f_t$} are smooth enough.}}
 \eqref{star_1} and \eqref{star_3} hold for  all $p>1$,  
 and if $h$ is supported in $V$,
  we get that  $t \mapsto \int h(x)\Theta(g(x)-a) \, d\rho_t$
is $\alpha$-H\"older for any 
$$\alpha <\frac{1}{p}
\bigl (1-\sup_t \frac{|\log \JJ_t|}{|\log \nu_{s,t}|}
\bigr )<1-\sup_t \frac{|\log \JJ_t|}{|\log \nu_{s,t}|}>0\, .$$
(Except if $|\det Df_t|\equiv 1$ for all $t$, this is weaker than  the claim in
\cite[Remark 2.3]{BKL}.)

 Note that $\JJ_t\in [(\nu_{u,t})^{d_u} (\bar \nu_{s,t})^{d_s},1]$,
 with  $\bar \nu_{s,t}$ defined as $\nu_{s,t}$,  replacing the sup over
 $v\in E^s\setminus \{0\}$ by an inf. Thus,
 $\sup_t \frac{|\log \JJ_t|}{|\log \nu_{s,t}|}\le \frac{|d_u \log \nu_u+d_s \log \bar \nu_s|}{|\log \nu_s|}<1$
 if 
 \begin{equation}
 \label{star_4} d_u \log \nu_u> d_s |\log \bar \nu_s|- |\log \nu_s|\, ,
 \end{equation} where
 $\nu_u=\inf_t \nu_{u,t}>1$ is the weakest expansion,
 $\nu_s=\sup_t \nu_{s,t}<1$ is the weakest contraction, and
 $\bar \nu_s=\inf_t \bar \nu_{s,t}\le \nu_s$ is the strongest contraction.
 
 We do not claim that the sufficient conditions \eqref{star_4} or   \eqref{star_3}
 are necessary, but some assumption is required to get fractional response: Stefano
 Galatolo pointed
 out at a conference in Warwick in April 2017 that, for a well-chosen perturbation $f_t$ of the standard
solenoidal attractor $f_0$, the map $t\mapsto \int_M \theta\,  d\rho_t$
is discontinuous at $t=0$ for certain  $\theta \in H^r_p$, since the SRB measure has atoms
in the stable direction.
(In this example, $d_u=1$ and $d_s=2$, with
$\JJ_0=1/2$ and $\nu_{s,0}=1/2$, so that
\eqref{star_3}
is violated.)

Changes in the proof of \cite[Theorem 2.1]{BKL} are limited to the following: 

$\bullet$Each $H^{\sigma}_q(M)$, for $1<q<\infty$
 and  $\sigma\in \real$   must be replaced by
$H^{\sigma}_q(\bar V)$.

$\bullet$6 lines below (2.3), p. 1210: ``for all real numbers $1<p'<\infty$ and $u-2<s<0<u<2$'' must be replaced
by ``for all real numbers $1<p'<\infty$ and all $u-2<s<0<u<2$  so
that \eqref{star_0} holds.''

$\bullet$3 lines above (2.4), p. 1210: ``enjoys the Perron-Frobenius spectral properties described above'' must be replaced by ``enjoys the Perron-Frobenius spectral properties described above,  if
\eqref{star_0}  holds when replacing $s$ by $s-1$ and $u$ by $u-1$.''

$\bullet$Lines 1-2 of p. 1212: ``(taking $s'<s<0$ close enough to zero)'' must be replaced by
``(taking $-r<s'<s<0$).''

$\bullet$2 lines after (2.8), p 1212: ``Fix   $|s|<r$
small'' must be replaced by ``Fix   $s_0<|s|<r$ close to 
$s_0:=\sup_t \frac{|\log \JJ_t|}{p|\log \nu_{s,t}|}<1$.''

$\bullet$Line after (2.11), p. 1212:  ``For each $0<\alpha<r$''
must be replaced  by ``For each $\alpha$ satisfying \eqref{star_2}''.  Last line of p. 1212, ``small enough   $\tilde r>|s|>0$'' must be replaced by ``$\tilde r>|s|>s_0$ close enough to $s_0$.''